\documentclass[11pt]{article}
\usepackage[left=1in,right=1in,top=1in,bottom=1in]{geometry}
\usepackage{times}
\usepackage{expl3}
\usepackage{cite}
\usepackage[table]{xcolor}
\usepackage{multirow}
\usepackage{stackengine} 
\usepackage{hhline}
\usepackage{lipsum}
\usepackage{titlesec}
\usepackage{wrapfig}
\usepackage{epsfig}
\usepackage{graphicx}
\usepackage{amsmath}
\usepackage[title]{appendix}
\usepackage{amssymb}
\usepackage{epstopdf}
\usepackage{boldline}
\usepackage{calligra}
\usepackage{url}
\usepackage{blindtext}

\newcommand{\define}{\stackrel{\mbox{\tiny def}}{=}}

\newtheorem{definition}{Definition}
\newtheorem{theorem}{Theorem}
\newtheorem{corollary}{Corollary}
\newtheorem{lemma}{Lemma}

\usepackage{mathtools}
\usepackage{epstopdf}
\usepackage{balance}
\usepackage{thmtools}
\usepackage{thm-restate}
\usepackage{hyperref}
\usepackage{cleveref}

\usepackage[ruled,vlined]{algorithm2e}
\include{pythonlisting}
\newcommand{\ostar}{\mathbin{\mathpalette\make@circled\star}}

\makeatletter
\newcommand{\removelatexerror}{\let\@latex@error\@gobble}
\makeatother
\makeatletter
\newcommand*{\rom}[1]{\expandafter\@slowromancap\romannumeral #1@}
\makeatother

\ExplSyntaxOn
\newcommand\latinabbrev[1]{
  \peek_meaning:NTF . {
    #1\@}%
  { \peek_catcode:NTF a {
      #1.\@ }%
    {#1.\@}}}
\ExplSyntaxOff

\titleclass{\subsubsubsection}{straight}[\subsubsection]

\begin{document}
\vspace{1cm}
\title{Random Double Tensors Integrals}\vspace{1.8cm}
\author{Shih~Yu~Chang 
\thanks{Shih Yu Chang is with the Department of Applied Data Science,
San Jose State University, San Jose, CA, U. S. A. (e-mail: {\tt
shihyu.chang@sjsu.edu}).
           }}

\maketitle

\begin{abstract}
In this work, we try to build a theory for random double tensor integrals (DTI). We begin with the definition of DTI and discuss how randomness structure is built upon DTI. Then, the tail bound of the unitarily invariant norm for the random DTI is established and this bound can help us to derive tail bounds of the unitarily invariant norm for various types of two tensors means, e.g., arithmetic mean, geometric mean, harmonic mean, and general mean. By associating DTI with perturbation formula, i.e., a formula to relate the tensor-valued function difference with respect the difference of the function input tensors, the tail bounds of the unitarily invariant norm for the Lipschitz estimate of tensor-valued function with random tensors as arguments are derived for vanilla case and quasi-commutator case, respectively.  We also establish the continuity property for random DTI in the sense of convergence in the random tensor mean, and we apply this continuity property to obtain the tail bound of the unitarily invariant norm for the derivative of the tensor-valued function.
\end{abstract}


\section{Introduction}\label{sec:Introduction}

In recent years, tensors have been applied to different applications in science and engineering~\cite{qi2017tensor, sidiropoulos2017tensor}. However, most of these applications assume that systems modelled by tensors are deterministic and such assumption is not always true and practical in problems involving tensor formulations. In recent years, more research results have pioneered some theories about random tensors~\cite{MR3616422,MR3783911,MR4140540}. One important question in random tensors is about concentration behavior of random tensors. In~\cite{chang2021convenient}, we extend Lapalace transform method and Lieb’s concavity theorem from matrices to tensors, and apply these tools to generalize the classical bounds associated with the names Chernoff, Bennett, and Bernstein from the scalar to the tensor setting. In~\cite{chang2021general}, this work extends previous work by considering the tail behavior of the top $k$-largest singular values of a function of the tensors summation, instead of the largest/smallest singular value of the tensors summation directly (identity function) explored in ~\cite{chang2021convenient}. Majorization and antisymmetric tensor product tools are main techniques utilized to establish inequalities for unitarily norms of multivariate tensors. Random tensors summation form discussed in~\cite{chang2021convenient,chang2021general} is linear form, i.e., each summand of random tensors with degree one. In works~\cite{chang2021hanson,chang2022hansongeneralized}, we extend the Hanson-Wright inequality for the maximum eigenvalue of the quadratic form of random Hermitian tensors under Einstein product. We separate the quadratic form of random tensors into diagonal summation and coupling (non-diagonal) summation parts. For the diagonal part, we can apply Bernstein inequality to bound the tail probability of the maximum eigenvalue of the summation of independent random Hermitian tensors directly. For coupling summation part, we have to apply decoupling method first, i.e., decoupling inequality to bound expressions with dependent random Hermitian tensors with independent random Hermitian tensors, before applying Bernstein inequality again to bound the tail probability of the maximum eigenvalue of the coupling summation of independent random Hermitian tensors. Previous works are based on tensors with Einstein products. Since Kilmer et al. introduced the new multiplication method between two third-order tensors around 2008 and third-order tensors with such multiplication structure are also called as T-product tensors~\cite{kilmer2008third}, T-product tensors have been applied to many fields in science and engineering, such as low-rank tensor approximation, signal processing, image feature extraction, machine learning, computer vision, and the multi-view clustering problem, etc. The discussion about concentration behaviors based on T-product tensors can also be found in~\cite{chang2022tP_I,chang2022tP_II}.

Inspired by operator mean theory (also called Kubo–Ando theory), we try to consider other operations besides $+$ (arithmetic mean) among tensors~\cite{hiai2010matrix}. The matrix mean for double operators can be expressed  by Eq.~(5:1:2) in~\cite{hiai2010matrix}, which has the same formation of double operator integral theory discussed in~\cite{skripka2019multilinear}. In this work, we begin to define double tensor integrals (DTI) and consider the tail bound for the unitarily invariant norm of random DTI, see Theorem~\ref{thm:unitary inv norm tail bound}. This bound can help us to establish tail bounds for various types of tensor means besides arithmetic mean. Since DTI can be used to express perturbation formula, i.e., an formula to relate the tensor-valued function difference with respect the difference of the function input tensors,  we establish Lipschitz estimate for random tensors by the tail bound format, see Theorem~\ref{thm:TB Lipschitz Estimates}. We also generalize Lipschitz estimate for random tensors with another quasi-commutator tensor, $\mathcal{D}$ by providing the tail bound for the unitarily invariant norm  of $\mathcal{D} \star_N f (\mathcal{A}) -  f (\mathcal{B}) \star_N \mathcal{D}$ in Theorem~\ref{thm:TB Lipschitz Estimates quasi}, where $\mathcal{A}, \mathcal{B}$ are random Hermitian tensors. We also establish a continuity for random DTI in the sense of convergence in the random tensor mean. This continuity property helps us to obtain the tail bounds for the unitarily invariant norm of the derivative of the tensor-valued function  under vanilla case and quasi-commutator case.

We define double tensor integrals (DTI) and randomness of DTI in Section~\ref{sec:Random Tensor Integral}. The tail bound for the unitarily invariant norm of random DTI and its applications to obtain various tail bounds for different types of double tensors means like arithmetic mean, geometric mean, harmonic mean, and general mean, are discussed in Section~\ref{sec:Tail Bound for Random Tensor Integral Norms}. In Section~\ref{sec:Tail Bound for Random Lipschitz Estimates}, we establish Lipschitz estimates for random tensors by the tail bound format for vanilla case and quasi-commutator case. We will establish continuity of random DTI based on the convergence in tensor mean of random Hermitian tensors in Section~\ref{sec:Continuity of Random Tensor Integral}. The application of DTI theory to acquire the tail bound for the unitarily invariant norm of the derivative of the tensor-valued function is presented by Section~\ref{sec:Applications of Tensor Integral}. Finally, the conclusions are given in Section~\ref{sec:Conclusions}.

\section{Random Double Tensor Integrals}\label{sec:Random Tensor Integral}

The purpose of this section is to define random double tensor integrals (DTI). We begin with the definition of DTI in Section~\ref{sec:Double Tensor Integrals}. In Section~\ref{sec:What Randomness Are Considered}, we will present what are randomness objects at DTI discussed at this work. 

Without loss of generality, one can partition the dimensions of a tensor into two groups, say $M$ and $N$ dimensions, separately. Thus, for two order-($M$+$N$) tensors: $\mathcal{X} \define (x_{i_1, \cdots, i_M, j_1, \cdots,j_N}) \in \mathbb{C}^{I_1 \times \cdots \times I_M\times
J_1 \times \cdots \times J_N}$ and $\mathcal{Y} \define (y_{i_1, \cdots, i_M, j_1, \cdots,j_N}) \in \mathbb{C}^{I_1 \times \cdots \times I_M\times
J_1 \times \cdots \times J_N}$, according to~\cite{MR3913666}, the \emph{tensor addition} $\mathcal{X} + \mathcal{Y}\in \mathbb{C}^{I_1 \times \cdots \times I_M\times
J_1 \times \cdots \times J_N}$ is given by 
\begin{eqnarray}\label{eq: tensor addition definition}
(\mathcal{X} + \mathcal{Y} )_{i_1, \cdots, i_M, j_1 , \cdots , j_N} &\define&
x_{i_1, \cdots, i_M, j_1 ,\cdots , j_N} \nonumber \\
& &+ y_{i_1, \cdots, i_M, j_1 ,\cdots , j_N}. 
\end{eqnarray}
On the other hand, for tensors $\mathcal{X} \define (x_{i_1, \cdots, i_M, j_1, \cdots,j_N}) \in \mathbb{C}^{I_1 \times \cdots \times I_M\times
J_1 \times \cdots \times J_N}$ and $\mathcal{Y} \define (y_{j_1, \cdots, j_N, k_1, \cdots,k_L}) \in \mathbb{C}^{J_1 \times \cdots \times J_N\times K_1 \times \cdots \times K_L}$, according to~\cite{MR3913666}, the \emph{Einstein product} (or simply referred to as \emph{tensor product} in this work) $\mathcal{X} \star_{N} \mathcal{Y} \in  \mathbb{C}^{I_1 \times \cdots \times I_M\times
K_1 \times \cdots \times K_L}$ is given by 
\begin{eqnarray}\label{eq: Einstein product definition}
\lefteqn{(\mathcal{X} \star_{N} \mathcal{Y} )_{i_1, \cdots, i_M,k_1, \cdots , k_L} \define} \nonumber \\ &&\sum\limits_{j_1, \cdots, j_N} x_{i_1, \cdots, i_M, j_1, \cdots,j_N}y_{j_1, \cdots, j_N, k_1, \cdots,k_L}. 
\end{eqnarray}

One can find more preliminary facts about tensors based on Einstein product in~\cite{chang2021convenient,MR3913666}. In the remaining of this paper, we will represent the scalar value $I_1 \times \cdots \times I_N$ by $\mathbb{I}_{1}^{N}$.

\subsection{Double Tensor Integrals}\label{sec:Double Tensor Integrals}

From Theorem 3.2 in~\cite{liang2019further}, every Hermitian tensor $\mathcal{H} \in  \mathbb{C}^{I_1 \times \cdots \times I_N \times I_1 \times \cdots \times I_N}$ has the following decomposition
\begin{eqnarray}\label{eq:Hermitian Eigen Decom}
\mathcal{H} &=& \sum\limits_{i=1}^{\mathbb{I}_{1}^{N}} \lambda_i \mathcal{U}_i \star_1 \mathcal{U}^{H}_i  \mbox{
~with~~$\langle \mathcal{U}_i, \mathcal{U}_i \rangle =1$ and $\langle \mathcal{U}_i, \mathcal{U}_j \rangle = 0$ for $i \neq j$,} \nonumber \\
&\define& \sum\limits_{i=1}^{\mathbb{I}_{1}^{N}} \lambda_i \mathcal{P}_{\mathcal{U}_i}
\end{eqnarray}
where $ \mathcal{U}_i \in  \mathbb{C}^{I_1 \times \cdots \times I_N \times 1}$, and the tensor $\mathcal{P}_{\mathcal{U}_i}$ is defined as $\mathcal{U}_i \star_1 \mathcal{U}^{H}_i$. The values $\lambda_i$ are named as \emph{eigevalues}. A Hermitian tensor with the decomposition shown by Eq.~\eqref{eq:Hermitian Eigen Decom} is named as \emph{eigen-decomposition}. A Hermitian tensor $\mathcal{H}$ is a positive definite (or positive semi-definite) tensor if all its eigenvalues are positive (or nonnegative).

Let $\mathcal{A}, \mathcal{B} \in  \mathbb{C}^{I_1 \times \cdots \times I_N \times I_1 \times \cdots \times I_N} $ be Hermitian tensors with the following eigen-decompositions:
\begin{eqnarray}
\mathcal{A} &=& \sum\limits_{i=1}^{\mathbb{I}_{1}^{N}} \lambda_i \mathcal{U}_i \star_1 \mathcal{U}^{H}_i
\define \sum\limits_{i=1}^{\mathbb{I}_{1}^{N}} \lambda_i \mathcal{P}_{\mathcal{U}_i},
\end{eqnarray}
and
\begin{eqnarray}
\mathcal{B} &=& \sum\limits_{j=1}^{\mathbb{I}_{1}^{N}} \mu_j \mathcal{V}_i \star_1 \mathcal{V}^{H}_j
\define \sum\limits_{j=1}^{\mathbb{I}_{1}^{N}} \mu_j  \mathcal{P}_{\mathcal{V}_j}.
\end{eqnarray}

We define \emph{double tensor integrals} (DTI) with respect to tensors $\mathcal{A}, \mathcal{B}$ and the function $\psi: \mathbb{R}^2 \rightarrow \mathbb{R}$, denoted as $T_{\mathcal{A}, \mathcal{B}, \psi}(\mathcal{X})$, which can be expressed as
\begin{eqnarray}\label{eq:Double Tensor Int Def}
T_{\mathcal{A}, \mathcal{B}, \psi}(\mathcal{X}) &=& \sum\limits_{i=1}^{\mathbb{I}_{1}^{N}} \sum\limits_{j=1}^{\mathbb{I}_{1}^{N}} \psi( \lambda_i, \mu_j)  \mathcal{P}_{\mathcal{U}_i} \star_N \mathcal{X} \star_N \mathcal{P}_{\mathcal{V}_j}.
\end{eqnarray}

\begin{lemma}\label{lma:tensor int algebraic properties}
Let $\psi, \phi: \mathbb{R}^2 \rightarrow \mathbb{R}$ be two functions, we have following relationships about $T_{\mathcal{A}, \mathcal{B}, \psi}(\mathcal{X})$:
Given $\psi$ is a constant function to one, we have 
\begin{eqnarray}\label{eq1:lma:tensor int algebraic properties}
T_{\mathcal{A}, \mathcal{B},  1} = \mathcal{I}.
\end{eqnarray}
We also have:
\begin{eqnarray}\label{eq2:lma:tensor int algebraic properties}
T_{\mathcal{A}, \mathcal{B},  \psi \phi }(\mathcal{X}) =  T_{\mathcal{A}, \mathcal{B}, \psi }(\mathcal{X}) \circ T_{\mathcal{A}, \mathcal{B}, \phi }(\mathcal{X}),
\end{eqnarray}
where $\circ$ is the  entrywise product (Hadamard product). Finally, we have
\begin{eqnarray}\label{eq3:lma:tensor int algebraic properties}
T_{\mathcal{A}, \mathcal{B}, a \psi + b \phi }(\mathcal{X}) =  aT_{\mathcal{A}, \mathcal{B}, \psi }(\mathcal{X})+b T_{\mathcal{A}, \mathcal{B}, \phi }(\mathcal{X}),  
\end{eqnarray}
where $a,b$ are two complex numbers. 
\end{lemma}
\textbf{Proof:}

Given a tensor $\mathcal{X} \in  \mathbb{C}^{I_1 \times \cdots \times I_N \times I_1 \times \cdots \times I_N}$ with orthgonal unitary tensors $\mathcal{U}_i$ and  $\mathcal{V}_j$ such that $\mathcal{A} = \sum\limits_{i=1}^{\mathbb{I}_{1}^{N}} \lambda_i \mathcal{U}_i \star_1 \mathcal{U}^{H}_i$  and $\mathcal{B} = \sum\limits_{j=1}^{\mathbb{I}_{1}^{N}} \mu_j \mathcal{U}_j  \star_1 \mathcal{U}^{H}_j$, we define
the scaler $x_{i, j}$ assocoated to $\mathcal{X}$ as 
\begin{eqnarray}
x_{i, j} &=& \langle \mathcal{X} \star_N \mathcal{V}_j, \mathcal{U}_i\rangle.
\end{eqnarray}
After selecting two specific indices $i'$ and $j'$, we have
\begin{eqnarray}\label{eq:matrix form}
\langle T_{\mathcal{A}, \mathcal{B}, \psi}(\mathcal{X})\star_N  \mathcal{V}_{j'}, \mathcal{U}_{i'}   \rangle
&=& \psi( \lambda_{i'}, \mu_{j'}) \langle \mathcal{X} \star_N \mathcal{V}_{j'}, \mathcal{U}_{i'} \rangle = \psi( \lambda_{i'}, \mu_{j'}) x_{i', j'}.
\end{eqnarray}
From Eq.~\eqref{eq:matrix form} and the Hadamard product properties, we have properties provided by Eqs.~\eqref{eq1:lma:tensor int algebraic properties},~\eqref{eq2:lma:tensor int algebraic properties} and~\eqref{eq3:lma:tensor int algebraic properties}.
$\hfill \Box$

\subsection{Random DTI}\label{sec:What Randomness Are Considered}

According to the DTI definition provided by Eq.~\eqref{eq:Double Tensor Int Def}, the random DTI considered in this work is to assume that tensors $\mathcal{A}, \mathcal{B}$ are random Hermitian tensors and the remaining parameters $\psi$ and $\mathcal{X}$ are deterministics. Therefore, we have the randomness at the following terms in Eq.~\eqref{eq:Double Tensor Int Def}:  $\psi( \lambda_i, \mu_j), \mathcal{P}_{\mathcal{U}_i}$ and $\mathcal{P}_{\mathcal{V}_j}$. If we are provided more detailed probability density functions for entries of random Hermitian tensors $\mathcal{A}, \mathcal{B}$, all bounds derived in this work can be improved with more dedicated expressions associated with parameters of probability density functions. 

\section{Tail Bound for Random Tensor Integral Norms}\label{sec:Tail Bound for Random Tensor Integral Norms}

\subsection{Unitarily Invariant Tensor Norms}\label{sec:Unitarily Invariant Tensor Norms}

Let us represent the Hermitian eigenvalues of a Hermitian tensor $\mathcal{H} \in \mathbb{C}^{I_1 \times \cdots \times I_N \times I_1 \times \cdots \times I_N} $ in decreasing order by the vector $\vec{\lambda}(\mathcal{H}) = (\lambda_1(\mathcal{H}), \cdots, \lambda_{\mathbb{I}_{1}^{N}}(\mathcal{H}))$. We use $\mathbb{R}_{\geq 0} (\mathbb{R}_{> 0})$ to represent a set of nonnegative (positive) real numbers. Let $\left\Vert \cdot \right\Vert_{\rho}$ be a unitarily invariant tensor norm, i.e., $\left\Vert \mathcal{H}\star_N \mathcal{U}\right\Vert_{\rho} = \left\Vert \mathcal{U}\star_N \mathcal{H}\right\Vert_{\rho} = \left\Vert \mathcal{H}\right\Vert_{\rho} $,  where $\mathcal{U}$ is any unitary tensor. Let $\rho : \mathbb{R}_{\geq 0}^{\mathbb{I}_{1}^{N}} \rightarrow \mathbb{R}_{\geq 0}$ be the corresponding gauge function that satisfies H$\ddot{o}$lder’s inequality so that 
\begin{eqnarray}\label{eq:def gauge func and general unitarily invariant norm}
\left\Vert \mathcal{H} \right\Vert_{\rho} = \left\Vert |\mathcal{H}| \right\Vert_{\rho} = \rho(\vec{\lambda}( | \mathcal{H} | ) ),
\end{eqnarray}
where $ |\mathcal{H}|  \define \sqrt{\mathcal{H}^H \star_N \mathcal{H}} $. 

We will provide several popular tensor norm examples which can be treated as special cases of unitarily invariant tensor norm.  The first one is Schatten $p$-norm for tensors, denoted as $\left\Vert \mathcal{X} \right\Vert_{p}$, is defined as:
\begin{eqnarray}\label{eq: Schatten p norm for tensors}
\left\Vert \mathcal{X}\right\Vert_{p} \define (\mathrm{Tr}|\mathcal{X}|^p )^{\frac{1}{p}},
\end{eqnarray}
where $ p \geq 1$. If $p=1$, it is the trace norm. 

The second one is $k$-trace norm, denoted as $\mathrm{Tr}_k[\mathcal{X}]$, defined by ~\cite{huang2020generalizing}. It is 
\begin{eqnarray}\label{eq: de k-trace norm for tensors}
\mathrm{Tr}_k[\mathcal{X}] \define \sum\limits_{1 \leq i_1 < i_2 < \cdots i_k \leq r} \lambda_{i_1} \lambda_{i_1}  \cdots \lambda_{i_k} 
\end{eqnarray}
where $ 1 \leq k \leq r$. If $k=1$, $\mathrm{Tr}_k[\mathcal{X}]$ is reduced as trace norm. 

The third one is Ky Fan like $k$-norm~\cite{fan1955some} for tensors. For $k \in \{1,2,\cdots, \mathbb{I}_1^{N} \}$, the Ky Fan $k$-norm~\cite{fan1955some} for tensors  $\mathcal{X} \in \mathbb{C}^{I_1 \times \cdots \times I_N \times I_1 \times \cdots \times I_N} $, denoted as $\left\Vert \mathcal{X}\right\Vert_{(k)}$, is defined as:
\begin{eqnarray}\label{eq: Ky Fan k norm for tensors}
\left\Vert \mathcal{X}\right\Vert_{(k)} \define \sum\limits_{i=1}^{k} \lambda_i(  |\mathcal{X}|  ).
\end{eqnarray}
If $k=1$,  the Ky Fan $k$-norm for tensors is the tensor operator norm, denoted as $ \left\Vert \mathcal{X} \right\Vert$. In this work, we will apply the symbol $\left\Vert \mathcal{X} \right\Vert_{\rho}$ to represent any unitarily invariant tensor norm for the tensor $\mathcal{X}$.

In the following theorem, we will present the tail bound of unitarily invariant tensor norm for a given tensor integral $T_{\mathcal{A}, \mathcal{B}, \psi}(\mathcal{X})$.  

\begin{theorem}\label{thm:unitary inv norm tail bound}
Let $\mathcal{A}, \mathcal{B} \in \mathbb{C}^{I_1 \times \cdots \times I_N \times I_1 \times \cdots \times I_N}$ be independent random Hermitian tensors with $\mathcal{A} = \sum\limits_{i=1}^{\mathbb{I}_{1}^{N}} \lambda_i \mathcal{U}_i \star_1 \mathcal{U}^{H}_i$  and $\mathcal{B} = \sum\limits_{j=1}^{\mathbb{I}_{1}^{N}} \mu_j \mathcal{U}_j  \star_1 \mathcal{U}^{H}_j$, then, for any $\theta > 0$, we have
\begin{eqnarray}\label{eq1:thm:unitary inv norm tail bound}
\mathrm{Pr}\left(   \left\Vert T_{\mathcal{A}, \mathcal{B}, \psi}(\mathcal{X})  \right\Vert_{\rho}  \geq \theta \right) \leq \frac{  \left( \mathbb{I}_1^N \right)^2 \left\Vert \mathcal{X} \right\Vert_{\rho} }{\theta}  \sum\limits_{i=1}^{\mathbb{I}_{1}^{N}}  \sum\limits_{j=1}^{\mathbb{I}_{1}^{N}} \mathbb{E}\left( \left\vert  \psi(\lambda_i, \mu_j)\right\vert \right), 
\end{eqnarray} 
where $\mathbb{E}$ is the expectation. 
\end{theorem}
\textbf{Proof:}

Since we have the following norm estimation for $\left\Vert T_{\mathcal{A}, \mathcal{B}, \psi}(\mathcal{X})  \right\Vert_{\rho}$:
\begin{eqnarray}\label{eq2:thm:unitary inv norm tail bound}
\left\Vert T_{\mathcal{A}, \mathcal{B}, \psi}(\mathcal{X})  \right\Vert_{\rho} &=& \left\Vert 
\sum\limits_{i=1}^{\mathbb{I}_{1}^{N}} \sum\limits_{j=1}^{\mathbb{I}_{1}^{N}} \psi( \lambda_i, \mu_j)  \mathcal{P}_{\mathcal{U}_i} \star_N \mathcal{X} \star_N \mathcal{P}_{\mathcal{V}_j}  \right\Vert_{\rho} \nonumber \\
& \leq_1 &
\sum\limits_{i=1}^{\mathbb{I}_{1}^{N}} \sum\limits_{j=1}^{\mathbb{I}_{1}^{N}} 
 \left\vert \psi( \lambda_i, \mu_j)  \right\vert  \left\Vert   \mathcal{P}_{\mathcal{U}_i} \star_N \mathcal{X} \star_N \mathcal{P}_{\mathcal{V}_j} \right\Vert_{\rho} \nonumber \\
& =_2 & \sum\limits_{i=1}^{\mathbb{I}_{1}^{N}} \sum\limits_{j=1}^{\mathbb{I}_{1}^{N}} 
 \left\vert \psi( \lambda_i, \mu_j)  \right\vert   \left\Vert \mathcal{X} \right\Vert_{\rho},
\end{eqnarray}
where $\leq_1$ comes from triangle inequality of the unitarily invariant norm and $=_2$ comes from the definition of the unitarily invariant norm. 

Then, we have the following bound for $\mathrm{Pr}\left(   \left\Vert T_{\mathcal{A}, \mathcal{B}, \psi}(\mathcal{X})  \right\Vert_{\rho}  \geq \theta \right)$
\begin{eqnarray}
\mathrm{Pr}\left(   \left\Vert T_{\mathcal{A}, \mathcal{B}, \psi}(\mathcal{X})  \right\Vert_{\rho}  \geq \theta \right) &\leq_1&  \mathrm{Pr}\left(   \sum\limits_{i=1}^{\mathbb{I}_{1}^{N}} \sum\limits_{j=1}^{\mathbb{I}_{1}^{N}} 
 \left\vert \psi( \lambda_i, \mu_j)  \right\vert   \left\Vert \mathcal{X} \right\Vert_{\rho}   \geq \theta \right) \nonumber \\
& \leq &  \sum\limits_{i=1}^{\mathbb{I}_{1}^{N}} \sum\limits_{j=1}^{\mathbb{I}_{1}^{N}} \mathrm{Pr} \left(
\left\vert \psi\left(\lambda_i, \mu_j \right) \right\vert \geq \frac{\theta}{\left\Vert \mathcal{X}\right\Vert_{\rho} \left( \mathbb{I}_1^N \right)^2  } \right) \nonumber \\
& \leq_2 & \frac{   \left\Vert \mathcal{X}\right\Vert_{\rho} \left( \mathbb{I}_1^N \right)^2    }{\theta}  \sum\limits_{i=1}^{\mathbb{I}_{1}^{N}} \sum\limits_{j=1}^{\mathbb{I}_{1}^{N}} \mathbb{E}\left( \left\vert \psi\left(\lambda_i, \mu_j \right) \right\vert \right), 
\end{eqnarray}
where $\leq_1$ comes from the inequality obtained by Eq.~\eqref{eq2:thm:unitary inv norm tail bound}, and the $\leq_2$ is based on Markov inequality. 
$\hfill \Box$

We will consider several important examples of $\phi$, which will represent different tensor means. Given
two random Hermitian tensors $\mathcal{A} = \sum\limits_{i=1}^{\mathbb{I}_{1}^{N}} \lambda_i \mathcal{U}_i \star_1 \mathcal{U}^{H}_i $ and $\mathcal{B}  = \sum\limits_{j=1}^{\mathbb{I}_{1}^{N}} \mu_j \mathcal{U}_j  \star_1 \mathcal{U}^{H}_j$, we use $p_{\lambda_i}(~)$ and $p_{\mu_j}(~)$ to represent the probability density functions for eigenvalues $\lambda_i$ and $\mu_j$, respectively. From Sec. 2.2. in~\cite{liang2019further}, one can apply unfolding technique to convert a random Hermitian tensor into a Hermitian matrix. If we have Assumption 3.1 in~\cite{ordonez2008ordered}, we are able to obtain the $i$-th eigenvalue distribution of a positive definite Hermitian tensor, see Corollary 3.3 in~\cite{ordonez2008ordered}.

\begin{corollary}[arithmetic mean]\label{cor:arithmetic mean}
Under conidtions provided by Theorem~\ref{thm:unitary inv norm tail bound}, if the function $\psi$ has the following form:
\begin{eqnarray}\label{eq1:cor:arithmetic mean}
\psi(x,y) = \frac{x + y}{2},
\end{eqnarray}
we have 
\begin{eqnarray}\label{eq1:cor:arithmetic mean}
\mathrm{Pr}\left(   \left\Vert T_{\mathcal{A}, \mathcal{B}, \psi = \frac{x + y}{2}}(\mathcal{X})  \right\Vert_{\rho}  \geq \theta \right) \leq \frac{  \left( \mathbb{I}_1^N \right)^2 \left\Vert \mathcal{X} \right\Vert_{\rho} }{2\theta}  \sum\limits_{i=1}^{\mathbb{I}_{1}^{N}}  \sum\limits_{j=1}^{\mathbb{I}_{1}^{N}} \left[\mathbb{E}\left( \left\vert \lambda_i \right\vert \right) + \mathbb{E}\left( \left\vert \mu_j \right\vert \right) \right], 
\end{eqnarray} 
where we have
\begin{eqnarray}
\mathbb{E}\left( \left\vert \lambda_i \right\vert \right) &=& \int_0^{\infty} t\left[ p_{\lambda_i}(t) + p_{\lambda_i}(-t)   \right] dt,
\end{eqnarray}
and
\begin{eqnarray}
\mathbb{E}\left( \left\vert \mu_j \right\vert \right) &=& \int_0^{\infty} t\left[ p_{\mu_j}(t) + p_{\mu_j}(-t)   \right] dt.
\end{eqnarray}
\end{corollary}
\textbf{Proof:}
The key is to evaluate $\mathbb{E}\left( \left\vert  \psi(\lambda_i, \mu_j)\right\vert \right)$, we have 
\begin{eqnarray}
\mathbb{E}\left( \left\vert  \psi(\lambda_i, \mu_j) \right\vert \right) &=& 
\mathbb{E}\left( \left\vert  \frac{ \lambda_i + \mu_j }{2} \right\vert \right) \nonumber \\
& \leq & \frac{1}{2} \mathbb{E}\left( \left\vert  \lambda_i  \right\vert \right) 
+ \frac{1}{2} \mathbb{E}\left( \left\vert  \mu_j  \right\vert \right).
\end{eqnarray}
This corollary is proved by the following fact for a random variable $Y = \left\vert X \right\vert$:
\begin{eqnarray}
  F_Y ( y )=
    \begin{cases}
      F_X ( y ) - F_X( -y ) & y \geq 0,\\
      0 & y < 0,
    \end{cases}  
\end{eqnarray}
where $F_Y$ and $F_X$ are CDFs of random variables $Y$ and $X$. 
$\hfill \Box$

\begin{corollary}[geometric mean]\label{cor:geometric mean}
Under conidtions provided by Theorem~\ref{thm:unitary inv norm tail bound} with the assumption that $\mathcal{A}$ and $\mathcal{B}$ are random positive definite tensors, if the function $\psi$ has the following form:
\begin{eqnarray}\label{eq1:cor:geometric mean}
\psi(x,y) = \sqrt{xy},
\end{eqnarray}
we have 
\begin{eqnarray}\label{eq1:cor:geometric mean}
\mathrm{Pr}\left(   \left\Vert T_{\mathcal{A}, \mathcal{B}, \psi = \sqrt{xy}}(\mathcal{X})  \right\Vert_{\rho}  \geq \theta \right) \leq \frac{  \left( \mathbb{I}_1^N \right)^2 \left\Vert \mathcal{X} \right\Vert_{\rho} }{\theta}  \sum\limits_{i=1}^{\mathbb{I}_{1}^{N}}  \sum\limits_{j=1}^{\mathbb{I}_{1}^{N}} \left[\mathbb{E}\left( \sqrt{\lambda_i }  \right) \mathbb{E}\left( \sqrt{\mu_j } \right) \right], 
\end{eqnarray} 
where we have
\begin{eqnarray}
\mathbb{E}\left( \sqrt{\lambda_i }  \right) &=& \int_0^{\infty} \sqrt{t}  p_{\lambda_i}(t) dt,
\end{eqnarray}
and
\begin{eqnarray}
\mathbb{E}\left( \sqrt{ \mu_j }  \right) &=& \int_0^{\infty} \sqrt{t}  p_{\mu_j}(t) dt.
\end{eqnarray}
\end{corollary}
\textbf{Proof:}
To evaluate $\mathbb{E}\left( \sqrt{\lambda_i \mu_j} \right)$, we have 
\begin{eqnarray}
\mathbb{E}\left( \left\vert  \psi(\lambda_i, \mu_j) \right\vert \right) &=& 
\mathbb{E}\left( \sqrt{\lambda_i \mu_j} \right) \nonumber \\
& = & \mathbb{E} \left( \sqrt{ \lambda_i } \right)  \mathbb{E}\left( \sqrt{ \mu_j } \right)
\end{eqnarray}
$\hfill \Box$

\begin{corollary}[harmonic mean]\label{cor:harmonic mean}
Under conidtions provided by Theorem~\ref{thm:unitary inv norm tail bound} with the assumption that $\mathcal{A}$ and $\mathcal{B}$ are random positive definite tensors, if the function $\psi$ has the following form:
\begin{eqnarray}\label{eq1:cor:harmonic mean}
\psi(x,y) = \frac{2}{x^{-1} + y^{-1}},
\end{eqnarray}
we have 
\begin{eqnarray}\label{eq1:cor:harmonic mean}
\mathrm{Pr}\left(   \left\Vert T_{\mathcal{A}, \mathcal{B}, \psi = \frac{2}{x^{-1} + y^{-1}} }(\mathcal{X})  \right\Vert_{\rho}  \geq \theta \right) &\leq & \frac{  2 \left( \mathbb{I}_1^N \right)^2 \left\Vert \mathcal{X} \right\Vert_{\rho} }{\theta}  \sum\limits_{i=1}^{\mathbb{I}_{1}^{N}}  \sum\limits_{j=1}^{\mathbb{I}_{1}^{N}}  \int_0^{\infty} \mathrm{H}_{\lambda_i}(t) \mathrm{H}_{\mu_j}(t) dt     , 
\end{eqnarray} 
where we have
\begin{eqnarray}
\mathrm{H}_{\lambda_i}(t) &=&  \int_0^{\infty}    \frac{e^{-tw}}{w^2}p_{\lambda_i}\left( 1/w \right) dw,
\end{eqnarray}
and
\begin{eqnarray}
\mathrm{H}_{\mu_j}(t)  &=&  \int_0^{\infty}    \frac{e^{-tz}}{z^2}p_{\mu_j}\left( 1/z \right) dz. 
\end{eqnarray}
\end{corollary}
\textbf{Proof:}
Let $X$ and $Y$ are two positive random variables with distribution fnuctions $f_X(x)$ and $f_Y(y)$. Then, we have
\begin{eqnarray}
\mathbb{E}\left( \frac{2}{X^{-1} + Y^{-1}} \right) &=_1&  \mathbb{E}\left( \frac{2}{W + Z} \right) \nonumber \\
&=&  2 \mathbb{E}\left( \int_0^{\infty}   \exp \left( -t (W+ Z )\right) dt \right) \nonumber \\
&=&  2 \int_0^{\infty}\mathbb{E}\left( \exp\left( -t ( W+Z  )\right) \right)d t  \nonumber \\
&=_2&  2 \int_0^{\infty}\mathbb{E}\left( \exp\left( -t W \right) \right)  \mathbb{E}\left( \exp\left( -t Z \right) \right) d t
\end{eqnarray}
where we set $W = X^{-1}$ and $Z = Y^{-1}$ at $=_1$ and we use indepedent assumptions of random variables $W = X^{-1}$ and $Z = Y^{-1}$ at $=_2$. This corollary is proved since we have 
following distribution functions for random variables $W$ and $Z$, which will correspond to random variables $\lambda_i$ and $\mu_j$, expressed as
\begin{eqnarray}
f(w) = \frac{1}{w^2}p_{\lambda_i}\left( 1/w\right),
\end{eqnarray}
and
\begin{eqnarray}
f(z) = \frac{1}{z^2}p_{\mu_j}\left( 1/z\right).
\end{eqnarray}
$\hfill \Box$

We have to prepare a lemma about the expectation of ratio between two depedent random variables before presenting the next corollary. 
\begin{lemma}\label{lma:exp of ratio dep RVs}
Given two random variables $X$ and $Y$ such that $Y \neq 0$ always, we have 
\begin{eqnarray}\label{eq1:lma:exp of ratio dep RVs}
\mathbb{E}\left( \frac{X}{Y} \right) = \frac{ \mathbb{E}\left( X \right) }{  \mathbb{E}\left( Y \right)  }
+ \lim\limits_{\epsilon \rightarrow 0} \sum\limits_{i=1}^{\infty}\frac{ (-1)^{i}\left[\mathbb{E}(X) \mathbb{E}(\acute{Y}^i) + \mathbb{E}( \acute{X}   \acute{Y}^i    ) \right]  }{ \prod\limits_{j=0}^i ( \mathbb{E}(Y)+ j \epsilon) },
\end{eqnarray}
where we define the following random variables:
\begin{eqnarray}
 \acute{X}  &\define& X - \mathbb{E}(X), \nonumber \\ 
 \acute{Y}  &\define& Y - \mathbb{E}(Y).  \nonumber \\ 
\end{eqnarray}
\end{lemma}
\textbf{Proof:}

We have the following expression about $\mathbb{E}(\frac{X}{Y})$, it is 
\begin{eqnarray}\label{eq1:lma:exp of ratio dep RVs}
\mathbb{E}(\frac{X}{Y}) &=& \mathbb{E}\left(\frac{ \mathbb{E}(X) }{ \mathbb{E}(Y)  } \times \frac{1 + \frac{\acute{X}}{\mathbb{E}(X)} }{    1 +  \frac{\acute{Y}}{\mathbb{E}(Y)}                   } \right) \nonumber \\
&=& \frac{ \mathbb{E}(X) }{ \mathbb{E}(Y)  } \mathbb{E}\left( \left( 1 + \frac{\acute{X}}{\mathbb{E}(X)} \right)   \left( 1 + \frac{\acute{Y}}{\mathbb{E}(Y)} \right)^{-1}   \right) \nonumber \\
&=& \frac{ \mathbb{E}(X) }{ \mathbb{E}(Y)  } \mathbb{E}\left(  \left( 1 + \frac{\acute{Y}}{\mathbb{E}(Y)} \right)^{-1}   \right) + \frac{1}{\mathbb{E}(Y)} \mathbb{E}\left(\acute{X}   \left( 1 + \frac{\acute{Y}}{\mathbb{E}(Y)} \right)^{-1}  \right)
\end{eqnarray}

Given a real function $g$, we have the following approximation form: 
\begin{eqnarray}\label{eq:Tayler exp}
g(a + x) = \lim\limits_{\epsilon \rightarrow 0} \sum\limits_{i=0}^{\infty} \frac{ x^i \Delta_{\epsilon}^{i} g(a)   }{\epsilon^{i} i !},
\end{eqnarray}
where $\Delta_{\epsilon}^{i}$ is the finite difference operator of degree $i$ and the step size $\epsilon$ is defined as
\begin{eqnarray}
\Delta_{\epsilon}^{i}g(a) = \sum\limits_{j=0}^{i}(-1)^{j} {i \choose j} g(a + (i - j) \epsilon).
\end{eqnarray}

If we apply Eq.~\eqref{eq:Tayler exp} to the following function at $ \acute{Y} = 0$
\begin{eqnarray}
g(\acute{Y}) = \left( 1 + \frac{\acute{Y}}{\mathbb{E}(Y)} \right)^{-1},
\end{eqnarray}
we will get 
\begin{eqnarray}\label{eq:Tayler exp 2}
\left( 1 + \frac{\acute{Y}}{\mathbb{E}(Y)} \right)^{-1}
&=& \lim\limits_{\epsilon \rightarrow 0} \sum\limits_{i=0}^{\infty} \frac{ (-1)^i (\acute{Y})^i \mathbb{E}(Y)}{  \prod\limits_{j=0}^{i}( \mathbb{E}(Y) + j \epsilon )    }.
\end{eqnarray}
This lemma is proved by applying Eq~\eqref{eq:Tayler exp 2} to Eq.~\eqref{eq1:lma:exp of ratio dep RVs}.
$\hfill \Box$

Following corollary is about the unitarily invariant norm tail bound for the two tensors mean general format. Note that we have logarithmic mean if $\alpha \rightarrow 1 $ in Eq.~\eqref{eq1:cor:general mean}.

\begin{corollary}[general mean]\label{cor:general mean}
Under conidtions provided by Theorem~\ref{thm:unitary inv norm tail bound} with the assumption that $\mathcal{A}$ and $\mathcal{B}$ are random positive definite tensors, if the function $\psi$ has the following form:
\begin{eqnarray}\label{eq1:cor:general mean}
\psi(x,y) = \frac{\alpha - 1}{\alpha} \frac{ x^{\alpha } - y^{\alpha } }{x^{\alpha - 1} - y^{\alpha - 1}},
\end{eqnarray}
where $\alpha \in \mathbb{R}$ and $x \neq y$~\footnote{If $x=y$, this situation has measure zero.}. Then we have 
\begin{eqnarray}\label{eq1:cor:general mean}
\mathrm{Pr}\left(   \left\Vert T_{\mathcal{A}, \mathcal{B}, \psi =  \frac{\alpha - 1}{\alpha} \frac{ x^{\alpha } - y^{\alpha } }{x^{\alpha - 1} - y^{\alpha - 1}} }(\mathcal{X})  \right\Vert_{\rho}  \geq \theta \right) \leq 
\frac{  \left( \mathbb{I}_1^N \right)^2 \left\Vert \mathcal{X} \right\Vert_{\rho}( \alpha -1 ) }{\theta \alpha} \times ~~~~~~~~~~~~~~~~~~~~~~~~~~~~~~~~~~~~~~~~~~~~~~~ 
\nonumber \\
\sum\limits_{i=1}^{\mathbb{I}_{1}^{N}}  \sum\limits_{j=1}^{\mathbb{I}_{1}^{N}} \left\{ \frac{ \mathbb{E}\left( X_{\lambda_i, \mu_j} \right) }{  \mathbb{E}\left( Y_{\lambda_i, \mu_j} \right)  }
+ \lim\limits_{\epsilon \rightarrow 0} \sum\limits_{\i=1}^{\infty}\frac{ (-1)^{\i}\left[\mathbb{E}(X_{\lambda_i, \mu_j}) \mathbb{E}((\acute{Y}_{\lambda_i, \mu_j})^{\i}) + \mathbb{E}( \acute{X}_{\lambda_i, \mu_j}   (\acute{Y}_{\lambda_i, \mu_j})^{\i}    ) \right]  }{ \prod\limits_{\j=0}^{\i} ( \mathbb{E}(Y_{\lambda_i, \mu_j})+ \j \epsilon) } \right\}, 
\end{eqnarray} 
where we have random variables $X_{\lambda_i, \mu_j}$ and $Y_{\lambda_i, \mu_j}$ defined by 
\begin{eqnarray}
X_{\lambda_i, \mu_j} &=& \lambda^{\alpha}_i  - \mu^{\alpha}_j,
\end{eqnarray}
and
\begin{eqnarray}
Y_{\lambda_i, \mu_j} &=& \lambda^{\alpha - 1}_i  - \mu^{\alpha - 1}_j.
\end{eqnarray}
\end{corollary}
\textbf{Proof:}

Since $\left\vert \psi(\lambda_i, \mu_j) \right\vert $ is 
\begin{eqnarray}\label{eq2:cor:general mean}
\psi(\lambda_i, \mu_j) &=& \frac{\alpha - 1}{ \alpha} \frac{  \lambda^{\alpha}_i  - \mu^{\alpha}_j   }{ \lambda^{\alpha - 1}_i  - \mu^{\alpha - 1}_j } = \frac{\alpha - 1}{ \alpha} \frac{  X_{\lambda_i, \mu_j}  }{  Y_{\lambda_i, \mu_j} }, 
\end{eqnarray}
this corollary is proved by applying Lemma~\ref{lma:exp of ratio dep RVs} to the expectation of Eq.~\eqref{eq2:cor:general mean}.  
$\hfill \Box$

Note that each of the following terms $ \mathbb{E}\left( X_{\lambda_i, \mu_j} \right),  \mathbb{E}\left( Y_{\lambda_i, \mu_j} \right),  \mathbb{E}((\acute{Y}_{\lambda_i, \mu_j})^{\i})$ and $\mathbb{E}( \acute{X}_{\lambda_i, \mu_j}   (\acute{Y}_{\lambda_i, \mu_j})^{\i}    ) $ can be evaluated exactly since we know all density distributions $p_{\lambda_i}(~)$ and $p_{\mu_j}(~)$. 

\section{Tail Bounds for Random Lipschitz Estimates}\label{sec:Tail Bound for Random Lipschitz Estimates}

In this section, we will try to provide tail bounds for the Lipschitz estimate for the unitarily invariant norm for a given function, which is the main result of this section. We will begin with the perturbation lemma. The vanilla case is discussed in Section~\ref{sec:Vanilla Case}. The case about considering quasi-commutator is provided by Section~\ref{sec:Quasi-Commutators Case}

\subsection{Vanilla Case}\label{sec:Vanilla Case}

We will begin by providing a perturbation formula for DTI. 

\begin{lemma}\label{lma:pert formula}
Let $\mathcal{A}, \mathcal{B} \in \mathbb{C}^{I_1 \times \cdots \times I_N \times I_1 \times \cdots \times I_N}$ be Hermitian tensors with $\mathcal{A} = \sum\limits_{i=1}^{\mathbb{I}_{1}^{N}} \lambda_i \mathcal{U}_i \star_1 \mathcal{U}^{H}_i$  and $\mathcal{B} = \sum\limits_{j=1}^{\mathbb{I}_{1}^{N}} \mu_j \mathcal{U}_j  \star_1 \mathcal{U}^{H}_j$. Also, let $f: \mathbb{R} \rightarrow \mathbb{R}$ such that $f'(x)$ exists, we define the bivariate function $f^{[1]}$ as  
\begin{eqnarray}\label{f[1] def:lma:pert formula}
f^{[1]}(x, y) =
    \begin{cases}
      \frac{f(x) - f(y)}{x - y}, & x \neq y ;\\
      f'(x), & x=y. 
    \end{cases}  
\end{eqnarray}
Then, we have 
\begin{eqnarray}\label{eq1:lma:pert formula}
f(\mathcal{A}) - f(\mathcal{B}) = T_{\mathcal{A}, \mathcal{B}, f^{[1]}}\left( \mathcal{A} - \mathcal{B} \right).  
\end{eqnarray}
\end{lemma}
\textbf{Proof:}

Since we have 
\begin{eqnarray}
\mathcal{A} = \sum\limits_{i=1}^{\mathbb{I}_{1}^{N}} \lambda_i \mathcal{U}_i \star_1 \mathcal{U}^{H}_i,
\end{eqnarray}
and
\begin{eqnarray}
\mathcal{B} = \sum\limits_{j=1}^{\mathbb{I}_{1}^{N}} \mu_j \mathcal{U}_j  \star_1 \mathcal{U}^{H}_j,
\end{eqnarray}
then, we will obtain the following
\begin{eqnarray}
f(\mathcal{A}) = \sum\limits_{i=1}^{\mathbb{I}_{1}^{N}} f(\lambda_i )  \mathcal{U}_i \star_1 \mathcal{U}^{H}_i,
\end{eqnarray}
and
\begin{eqnarray}
f(\mathcal{B}) = \sum\limits_{j=1}^{\mathbb{I}_{1}^{N}} f(\mu_j ) \mathcal{U}_j  \star_1 \mathcal{U}^{H}_j. 
\end{eqnarray}

Because, we also have
\begin{eqnarray}
f(\mathcal{A})\star_N \mathcal{U}_i   =  f(\lambda_i )  \star_N \mathcal{U}_i, 
\end{eqnarray}
and
\begin{eqnarray}
f(\mathcal{B})\star_N \mathcal{V}_j   =  f(\mu_j)  \star_N \mathcal{V}_j, 
\end{eqnarray}
then, we have
\begin{eqnarray}\label{eq2:lma:pert formula}
\langle \left( f(\mathcal{A}) - f(\mathcal{B})  \star_N \mathcal{V}_j,  \mathcal{U}_i \right) \rangle
&=& f^{[1]}(\lambda_i, \mu_j) \langle \left( \mathcal{A} - \mathcal{B} \right) \star_N \mathcal{V}_j, \mathcal{U}_i  \rangle.
\end{eqnarray}

By applying Eq.~\eqref{eq:Double Tensor Int Def} definition and Eq.~\eqref{eq:matrix form} to Eq.~\eqref{eq2:lma:pert formula}, this Lemma is proved. 
$\hfill \Box$

From perturbation formula given by Lemma~\ref{lma:pert formula}, we can have 
the following theoem about the tail bounds of the unitarily invariant norm for the Lipschitz estimate of tensor-valued function with random tensors as inputs. 

\begin{theorem}\label{thm:TB Lipschitz Estimates}
Let $\mathcal{A}, \mathcal{B} \in \mathbb{C}^{I_1 \times \cdots \times I_N \times I_1 \times \cdots \times I_N}$ be independent random Hermitian tensors with $\mathcal{A} = \sum\limits_{i=1}^{\mathbb{I}_{1}^{N}} \lambda_i \mathcal{U}_i \star_1 \mathcal{U}^{H}_i$  and $\mathcal{B} = \sum\limits_{j=1}^{\mathbb{I}_{1}^{N}} \mu_j \mathcal{U}_j  \star_1 \mathcal{U}^{H}_j$, moreover, we are given a real valud funtion $f(x)$ for $x \in \mathbb{R}$ such that $\left\vert f^{[1]}(x, y) \right\vert$ is bounded by a positive number denoted as $\Omega_{f^{[1]}}$. Then, for any $\theta > 0$, we have
\begin{eqnarray}\label{eq1:thm:TB Lipschitz Estimates}
\mathrm{Pr}\left(  \left\Vert f(\mathcal{A}) -  f(\mathcal{B})  \right\Vert_{\rho}  \geq \theta \right) \leq \frac{  \left( \mathbb{I}_1^N \right)^2  \Omega_{f^{[1]}}  }{\theta} \mathbb{E}\left( \left\Vert \mathcal{A} - \mathcal{B} \right\Vert_{\rho} \right).
\end{eqnarray} 
\end{theorem}
\textbf{Proof:}

We have
\begin{eqnarray}\label{eq2:thm:TB Lipschitz Estimates}
\left\Vert f(\mathcal{A}) - f(\mathcal{B}) \right\Vert_{\rho}&=& \left\Vert T_{\mathcal{A}, \mathcal{B}, f^{[1]}}\left( \mathcal{A} - \mathcal{B} \right)  \right\Vert_{\rho} \nonumber \\
 &=&  \left\Vert \sum\limits_{i=1}^{\mathbb{I}_{1}^{N}} \sum\limits_{j=1}^{\mathbb{I}_{1}^{N}} f^{[1]}( \lambda_i, \mu_j)  \mathcal{P}_{\mathcal{U}_i} \star_N \left( \mathcal{A} - \mathcal{B} \right) \star_N \mathcal{P}_{\mathcal{V}_j}  \right\Vert_{\rho} \nonumber \\
& \leq_1 &
\sum\limits_{i=1}^{\mathbb{I}_{1}^{N}} \sum\limits_{j=1}^{\mathbb{I}_{1}^{N}} 
 \left\vert  f^{[1]}( \lambda_i, \mu_j)    \right\vert  \left\Vert   \mathcal{P}_{\mathcal{U}_i} \star_N \left( \mathcal{A} - \mathcal{B} \right) \star_N \mathcal{P}_{\mathcal{V}_j} \right\Vert_{\rho} \nonumber \\
& =_2 & \sum\limits_{i=1}^{\mathbb{I}_{1}^{N}} \sum\limits_{j=1}^{\mathbb{I}_{1}^{N}} 
 \left\vert f^{[1]}( \lambda_i, \mu_j)   \right\vert   \left\Vert \mathcal{A} - \mathcal{B}  \right\Vert_{\rho}  \nonumber \\ 
& \leq_3 & \left( \mathbb{I}_1^N \right)^2 \Omega_{f^{[1]}}   \left\Vert \mathcal{A} - \mathcal{B}  \right\Vert_{\rho},
\end{eqnarray}
where $\leq_1$ comes from triangle inequality of the unitarily invariant norm, $=_2$ comes from the definition of the unitarily invariant norm, and $\leq_3$ is due to that $\left\vert f^{[1]}( \lambda_i, \mu_j)   \right\vert \leq \Omega_{f^{[1]}}$ (mean value theorem of divide difference). 

Then, we have the following bound for $\mathrm{Pr}\left(  \left\Vert f(\mathcal{A}) - f(\mathcal{B}) \right\Vert_{\rho}  \geq \theta \right)$
\begin{eqnarray}
\mathrm{Pr}\left(   \left\Vert f(\mathcal{A}) - f(\mathcal{B}) \right\Vert_{\rho}  \geq \theta \right) &\leq_1&  \mathrm{Pr}\left(   \left( \mathbb{I}_1^N \right)^2 \Omega_{f^{[1]}}   \left\Vert \mathcal{A} - \mathcal{B} \right\Vert_{\rho}  \geq \theta \right) \nonumber \\
& \leq_2 & \frac{   \left( \mathbb{I}_1^N \right)^2 \Omega_{f^{[1]}}      }{\theta}  \mathbb{E}\left( \left\vert \mathcal{A} - \mathcal{B} \right\vert_{\rho} \right), 
\end{eqnarray}
where $\leq_1$ comes from the inequality obtained by Eq.~\eqref{eq2:thm:TB Lipschitz Estimates}, and the $\leq_2$ is based on Markov inequality. 
$\hfill \Box$

Since the upper bound $\Omega_{f^{[1]}}$ depends on $f(x)$, we will consider the following two corollaries about special types of the function $f(x)$. 

\begin{corollary}[Lipschitz estimate for polynomial fnuctions]\label{cor:TB Lipschitz Estimates:poly}
Let $\mathcal{A}, \mathcal{B} \in \mathbb{C}^{I_1 \times \cdots \times I_N \times I_1 \times \cdots \times I_N}$ be independent random Hermitian tensors with $\mathcal{A} = \sum\limits_{i=1}^{\mathbb{I}_{1}^{N}} \lambda_i \mathcal{U}_i \star_1 \mathcal{U}^{H}_i$  and $\mathcal{B} = \sum\limits_{j=1}^{\mathbb{I}_{1}^{N}} \mu_j \mathcal{U}_j  \star_1 \mathcal{U}^{H}_j$. Besides, we are given a real polynomial funtion $f(x)$ with degree $m$ over an interval $[a, b]$. For any $\theta > 0$, we have
\begin{eqnarray}\label{eq1:cor:TB Lipschitz Estimates:poly}
\mathrm{Pr}\left(  \left\Vert f(\mathcal{A}) -  f(\mathcal{B})  \right\Vert_{\rho}  \geq \theta \right) \leq \frac{  \left( \mathbb{I}_1^N \right)^2  \sum\limits_{k=1}^m \frac{ \left\vert f^{(k)}(x^{*}_{(k)})\right\vert(b-a)^k  }{k !}  }{\theta} \mathbb{E}\left( \left\Vert \mathcal{A} - \mathcal{B} \right\Vert_{\rho} \right),
\end{eqnarray} 
where $x^{*}_{(k)}$ is the maximizer to reach the maximum value for the function $\left\vert f^{(k)}(x)\right\vert$, i.e., the absolute value of the $k$-th derivative, in the interval $[a, b]$. 
\end{corollary}
\textbf{Proof:}

If we perform Tayler expansion for the function $f$ at $x$, we have
\begin{eqnarray}\label{eq1:cor:TB Lipschitz Estimates:poly}
f(y) = f(x) + f'(x)(y-x) + f''(x)\frac{(y-x)^2}{2!} +  f'''(x)\frac{(y-x)^3}{3!} + \ldots,
\end{eqnarray}
which is equivalent to have
\begin{eqnarray}
\frac{f(y) - f(x)}{y - x} =f'(x) + f''(x)\frac{(y-x)}{2!} +  f'''(x)\frac{(y-x)^2}{3!} + \ldots.
\end{eqnarray}

Because the polynomial function $f$ has degree $m$, from Eq.~\eqref{eq1:cor:TB Lipschitz Estimates:poly}, we have the following bound from triangle inequality:
\begin{eqnarray}
\left\vert \frac{f(y) - f(x)}{y - x} \right\vert &=& \left\vert   \sum\limits_{k=1}^{m} \frac{  f^{(k)}(x)(y-x)^k}{k!} \right\vert \nonumber \\
&\leq & \sum\limits_{k=1}^{m}  \frac{ \left\vert f^{(k)}(x^{*}_{(k)})\right\vert(b-a)^k  }{k !}, 
\end{eqnarray}
where $x^{*}_{(k)}$ is the maximizer to reach the maximum value for the function $\left\vert f^{(k)}(x^{*}_{(k)})\right\vert$ in the interval $[a, b]$. This corrollary is proved by Theorem~\ref{thm:TB Lipschitz Estimates}.
$\hfill \Box$

Following corollary is about Lipschitz estimate for polygamma functions. Recall that a digamma function $\omega(x)$
is defined as
\begin{eqnarray}\label{eq:polygamma func def}
\omega(x) = \frac{\Gamma'(x)}{\Gamma(x)},
\end{eqnarray}
where $\Gamma(x) = \int_0^{\infty} e^{-t} t^{x-1} dt$~\cite{qi2008new}. Then polygamma functions are defined as the $k$-th derivative $\omega^{(k)}(x)$ for any $k \in \mathbb{N}$.

\begin{corollary}[Lipschitz estimate for polygamma functions]\label{cor:TB Lipschitz Estimates:polygamma}
Let $\mathcal{A}, \mathcal{B} \in \mathbb{C}^{I_1 \times \cdots \times I_N \times I_1 \times \cdots \times I_N}$ be independent random positive definite tensors with $\mathcal{A} = \sum\limits_{i=1}^{\mathbb{I}_{1}^{N}} \lambda_i \mathcal{U}_i \star_1 \mathcal{U}^{H}_i$  and $\mathcal{B} = \sum\limits_{j=1}^{\mathbb{I}_{1}^{N}} \mu_j \mathcal{U}_j  \star_1 \mathcal{U}^{H}_j$. Besides, we are given a polygamma funtion $\omega^{(k)}(x)$  for any $k \in \mathbb{N}$ and $x > 0$.  For any $\theta > 0$, we have
\begin{eqnarray}\label{eq1:cor:TB Lipschitz Estimates:polygamma}
\mathrm{Pr}\left(  \left\Vert \omega^{(k)}(\mathcal{A}) -  \omega^{(k)}(\mathcal{B})  \right\Vert_{\rho}  \geq \theta \right) \leq \frac{  \left( \mathbb{I}_1^N \right)^2 \Omega_{\omega^{[k+1]}  }}{\theta} \mathbb{E}\left( \left\Vert \mathcal{A} - \mathcal{B} \right\Vert_{\rho} \right),
\end{eqnarray} 
where 
\begin{eqnarray}
\Omega_{\omega^{[k+1]}}  & = &  \max\left\{ \omega^{(k+1)}\left(\frac{1}{e}\left( \frac{(b^*)^{b^*}}{(a^*)^{a^*}}\right)^{\frac{1}{b^*-a^*}} \right)   ,  \omega^{(k+1)}(x^*) \right\}.
\end{eqnarray}
The values $a^*$ and $b^*$ are the maximizers of the function $\omega^{(k+1)}\left(\frac{1}{e}\left( \frac{b^{b}}{a^{a}}\right)^{\frac{1}{b-a}} \right)$ given $a \neq b$ and $a, b > 0$. The value $x^*$ is the maximizer of the function $\omega^{(k+1)}(x)$ given $x > 0$. 
\end{corollary}
\textbf{Proof:}

This corrollary is proved according to Theorem~\ref{thm:TB Lipschitz Estimates}  by applying Theorem 1 from~\cite{qi2008new} to Eq.~\eqref{eq:polygamma func def}. 
$\hfill \Box$

\subsection{Quasi-Commutator Case}\label{sec:Quasi-Commutators Case}

In this section, we will extend the perturbation formula provided by Lemma~\ref{lma:pert formula quasi} to the quasi-commutator $\mathcal{D} \star \mathcal{A} -  \mathcal{B} \star \mathcal{D}$. Tail bouns based on the quasi-commutator $\mathcal{D} \star \mathcal{A} -  \mathcal{B} \star \mathcal{D}$ will be given in this section.

\begin{lemma}\label{lma:pert formula quasi}
Let $\mathcal{A}, \mathcal{B} \in \mathbb{C}^{I_1 \times \cdots \times I_N \times I_1 \times \cdots \times I_N}$ be random Hermitian tensors with $\mathcal{A} = \sum\limits_{i=1}^{\mathbb{I}_{1}^{N}} \lambda_i \mathcal{U}_i \star_1 \mathcal{U}^{H}_i$  and $\mathcal{B} = \sum\limits_{j=1}^{\mathbb{I}_{1}^{N}} \mu_j \mathcal{V}_j  \star_1 \mathcal{V}^{H}_j$.  Besides, we have the function $f: \mathbb{R} \rightarrow \mathbb{R}$ such that $f'(x)$ exists, 
then, we have 
\begin{eqnarray}\label{eq1:lma:pert formula quasi}
\mathcal{D} \star_N f(\mathcal{A}) - f(\mathcal{B}) \star_N \mathcal{D} = T_{\mathcal{A}, \mathcal{B}, f^{[1]}}\left( \mathcal{D} \star_N \mathcal{A} -  \mathcal{B} \star_N \mathcal{D} \right).  
\end{eqnarray}
where $f^{[1]}$ has been defined by Eq.~\eqref{f[1] def:lma:pert formula}.
\end{lemma}
\textbf{Proof:}

Since we have 
\begin{eqnarray}
\mathcal{A} = \sum\limits_{i=1}^{\mathbb{I}_{1}^{N}} \lambda_i \mathcal{U}_i \star_1 \mathcal{U}^{H}_i,
\end{eqnarray}
and
\begin{eqnarray}
\mathcal{B} = \sum\limits_{j=1}^{\mathbb{I}_{1}^{N}} \mu_j \mathcal{V}_j  \star_1 \mathcal{V}^{H}_j,
\end{eqnarray}
then, we will obtain the following
\begin{eqnarray}
\mathcal{D} \star_N \mathcal{A} =  T_{\mathcal{A},\mathcal{I}, \lambda_i}  \star_N  \mathcal{D},
\end{eqnarray}
and
\begin{eqnarray}
\mathcal{B} \star_N \mathcal{D} =  T_{\mathcal{I}, \mathcal{B}, \mu_j}  \star_N  \mathcal{D}.
\end{eqnarray}

By applying Lemma~\ref{lma:tensor int algebraic properties}, we have
\begin{eqnarray}
T_{\mathcal{A}, \mathcal{B}, f^{[1]}}\left( \mathcal{D} \star_N \mathcal{A}  - \mathcal{B} \star_N \mathcal{D}\right) &=& T_{\mathcal{A}, \mathcal{B}, f^{[1]}}\left( T_{\mathcal{A},\mathcal{I}, \lambda_i}  \star_N  \mathcal{D}   - T_{\mathcal{I}, \mathcal{B}, \mu_j}  \star_N  \mathcal{D}\right) \nonumber \\
&=&  T_{\mathcal{A}, \mathcal{B}, f^{[1]}( \lambda_i , \mu_j)}\left(  \mathcal{D}\right)  \nonumber \\
&=&  T_{\mathcal{A}, \mathcal{B}, f(\lambda_i) - f(\mu_j)}\left(  \mathcal{D}\right) \nonumber \\
&=&  \mathcal{D} \star_N f (\mathcal{A}) -  f (\mathcal{B}) \star_N \mathcal{D}.
\end{eqnarray}
$\hfill \Box$

Following theorem is the tail bound for $\mathcal{D} \star_N f(\mathcal{A}) - f(\mathcal{B}) \star_N \mathcal{D}$.

\begin{theorem}\label{thm:TB Lipschitz Estimates quasi}
Let $\mathcal{A}, \mathcal{B} \in \mathbb{C}^{I_1 \times \cdots \times I_N \times I_1 \times \cdots \times I_N}$ be independent random Hermitian tensors with $\mathcal{A} = \sum\limits_{i=1}^{\mathbb{I}_{1}^{N}} \lambda_i \mathcal{U}_i \star_1 \mathcal{U}^{H}_i$  and $\mathcal{B} = \sum\limits_{j=1}^{\mathbb{I}_{1}^{N}} \mu_j \mathcal{V}_j  \star_1 \mathcal{V}^{H}_j$, moreover, we are given a real valud funtion $f(x)$ for $x \in \mathbb{R}$ such that $f'(x)$ exists and $\left\vert f'(x) \right\vert$ is bounded by a positive number denoted as $\Omega_{f^{[1]}}$. Then, for any $\theta > 0$, we have
\begin{eqnarray}\label{eq1:thm:TB Lipschitz Estimates}
\mathrm{Pr}\left(  \left\Vert \mathcal{D} \star_N f(\mathcal{A}) -  f(\mathcal{B}) \star_N \mathcal{D}  \right\Vert_{\rho}  \geq \theta \right) \leq \frac{  \left( \mathbb{I}_1^N \right)^2  \Omega_{f^{[1]}}  }{\theta} \mathbb{E}\left( \left\Vert \mathcal{D} \star_N \mathcal{A} - \mathcal{B} \star_N \mathcal{D} \right\Vert_{\rho} \right).
\end{eqnarray} 
\end{theorem}
\textbf{Proof:}

From Lemma~\ref{lma:pert formula quasi}, we have
\begin{eqnarray}\label{eq2:thm:TB Lipschitz Estimates quasi}
\left\Vert \mathcal{D} \star_N f(\mathcal{A}) - f(\mathcal{B}) \star_N \mathcal{D} \right\Vert_{\rho}&=& \left\Vert T_{\mathcal{A}, \mathcal{B}, f^{[1]}}\left( \mathcal{D} \star_N \mathcal{A} - \mathcal{B} \star_N \mathcal{D} \right)  \right\Vert_{\rho} \nonumber \\
 &=&  \left\Vert \sum\limits_{i=1}^{\mathbb{I}_{1}^{N}} \sum\limits_{j=1}^{\mathbb{I}_{1}^{N}} f^{[1]}( \lambda_i, \mu_j)  \mathcal{P}_{\mathcal{U}_i} \star_N \left(  \mathcal{D} \star_N \mathcal{A} - \mathcal{B} \star_N \mathcal{D} \right) \star_N \mathcal{P}_{\mathcal{V}_j}  \right\Vert_{\rho} \nonumber \\
& \leq_1 &
\sum\limits_{i=1}^{\mathbb{I}_{1}^{N}} \sum\limits_{j=1}^{\mathbb{I}_{1}^{N}} 
 \left\vert  f^{[1]}( \lambda_i, \mu_j)    \right\vert  \left\Vert   \mathcal{P}_{\mathcal{U}_i} \star_N \left(  \mathcal{D} \star_N \mathcal{A} - \mathcal{B} \star_N \mathcal{D} \right) \star_N \mathcal{P}_{\mathcal{V}_j} \right\Vert_{\rho} \nonumber \\
& =_2 & \sum\limits_{i=1}^{\mathbb{I}_{1}^{N}} \sum\limits_{j=1}^{\mathbb{I}_{1}^{N}} 
 \left\vert f^{[1]}( \lambda_i, \mu_j)   \right\vert   \left\Vert  \mathcal{D} \star_N \mathcal{A} - \mathcal{B} \star_N \mathcal{D} \right\Vert_{\rho}  \nonumber \\ 
& \leq_3 & \left( \mathbb{I}_1^N \right)^2 \Omega_{f^{[1]}}   \left\Vert  \mathcal{D} \star_N \mathcal{A} - \mathcal{B} \star_N \mathcal{D} \right\Vert_{\rho},
\end{eqnarray}
where $\leq_1$ comes from triangle inequality of the unitarily invariant norm, $=_2$ comes from the definition of the unitarily invariant norm, and $\leq_3$ is due to that $\left\vert f^{[1]}( \lambda_i, \mu_j)   \right\vert \leq \Omega_{f^{[1]}}$ (mean value theorem of divide difference). 

Then, we have the following bound for $\mathrm{Pr}\left(  \left\Vert \mathcal{D} \star_N f(\mathcal{A}) - f(\mathcal{B})  \star_N \mathcal{D}  \right\Vert_{\rho}  \geq \theta \right)$
\begin{eqnarray}
\mathrm{Pr}\left(   \left\Vert \mathcal{D} \star_N f(\mathcal{A}) - f(\mathcal{B}) \star_N \mathcal{D} \right\Vert_{\rho}  \geq \theta \right) &\leq_1&  \mathrm{Pr}\left(   \left( \mathbb{I}_1^N \right)^2 \Omega_{f^{[1]}}   \left\Vert  \mathcal{D} \star_N \mathcal{A} - \mathcal{B} \star_N \mathcal{D} \right\Vert_{\rho}  \geq \theta \right) \nonumber \\
& \leq_2 & \frac{   \left( \mathbb{I}_1^N \right)^2 \Omega_{f^{[1]}}      }{\theta}  \mathbb{E}\left( \left\Vert  \mathcal{D} \star_N \mathcal{A} - \mathcal{B} \star_N \mathcal{D}\right\Vert_{\rho} \right), 
\end{eqnarray}
where $\leq_1$ comes from the inequality obtained by Eq.~\eqref{eq2:thm:TB Lipschitz Estimates quasi}, and the $\leq_2$ is based on Markov inequality. 
$\hfill \Box$

Similar to Corollaries~\ref{cor:TB Lipschitz Estimates:poly} and~\ref{cor:TB Lipschitz Estimates:polygamma}, we have following two corollaries for Lipschitz estimate for polynomial and polygamma under quasi-commutator case.

\begin{corollary}[Lipschitz estimate for polynomial fnuctions, quasi-commutator case]\label{cor:TB Lipschitz Estimates:poly quasi}

Let $\mathcal{A}, \mathcal{B} \in $ \\ $ \mathbb{C}^{I_1 \times \cdots \times I_N \times I_1 \times \cdots \times I_N}$ be independent random Hermitian tensors with $\mathcal{A} = \sum\limits_{i=1}^{\mathbb{I}_{1}^{N}} \lambda_i \mathcal{U}_i \star_1 \mathcal{U}^{H}_i$  and $\mathcal{B} = \sum\limits_{j=1}^{\mathbb{I}_{1}^{N}} \mu_j \mathcal{V}_j  \star_1 \mathcal{V}^{H}_j$. Besides, we are given a real polynomial funtion $f(x)$ with degree $m$ over an interval $[a, b]$. For any $\theta > 0$, we have
\begin{eqnarray}\label{eq1:cor:TB Lipschitz Estimates:poly quasi}
\mathrm{Pr}\left(  \left\Vert \mathcal{D} \star_N f(\mathcal{A}) -  f(\mathcal{B}) \star_N  \mathcal{D}  \right\Vert_{\rho}  \geq \theta \right) \leq \frac{  \left( \mathbb{I}_1^N \right)^2  \sum\limits_{k=1}^m \frac{ \left\vert f^{(k)}(x^{*}_{(k)})\right\vert(b-a)^k  }{k !}  }{\theta} \mathbb{E}\left( \left\Vert \mathcal{D} \star_N \mathcal{A} - \mathcal{B} \star_N  \mathcal{D} \right\Vert_{\rho} \right),
\end{eqnarray} 
where $x^{*}_{(k)}$ is the maximizer to reach the maximum value for the function $\left\vert f^{(k)}(x^{*}_{(k)})\right\vert$ in the interval $[a, b]$. 
\end{corollary}
\textbf{Proof:}

If we perform Tayler expansion for the function $f$ at $x$, we have
\begin{eqnarray}\label{eq1:cor:TB Lipschitz Estimates:poly quasi}
f(y) = f(x) + f'(x)(y-x) + f''(x)\frac{(y-x)^2}{2!} +  f'''(x)\frac{(y-x)^3}{3!} + \ldots,
\end{eqnarray}
which is equivalent to have
\begin{eqnarray}
\frac{f(y) - f(x)}{y - x} =f'(x) + f''(x)\frac{(y-x)}{2!} +  f'''(x)\frac{(y-x)^2}{3!} + \ldots.
\end{eqnarray}

Because the polynomial function $f$ has degree $m$, from Eq.~\eqref{eq1:cor:TB Lipschitz Estimates:poly}, we have the following bound from triangle inequality:
\begin{eqnarray}
\left\vert \frac{f(y) - f(x)}{y - x} \right\vert &=& \left\vert   \sum\limits_{k=1}^{m} \frac{  f^{(k)}(x)(y-x)^k}{k!} \right\vert \nonumber \\
&\leq & \sum\limits_{k=1}^{m}  \frac{ \left\vert f^{(k)}(x^{*}_{(k)})\right\vert(b-a)^k  }{k !}, 
\end{eqnarray}
where $x^{*}_{(k)}$ is the maximizer to reach the maximum value for the function $\left\vert f^{(k)}(x^{*}_{(k)})\right\vert$ in the interval $[a, b]$. This corrollary is proved by Theorem~\ref{thm:TB Lipschitz Estimates quasi}.
$\hfill \Box$

\begin{corollary}[Lipschitz estimate for polygamma functions: quasi-commutator case]\label{cor:TB Lipschitz Estimates:polygamma quasi}

Let $\mathcal{A}, \mathcal{B} \in $ \\ $\mathbb{C}^{I_1 \times \cdots \times I_N \times I_1 \times \cdots \times I_N}$ be independent random positive definite tensors with $\mathcal{A} = \sum\limits_{i=1}^{\mathbb{I}_{1}^{N}} \lambda_i \mathcal{U}_i \star_1 \mathcal{U}^{H}_i$  and $\mathcal{B} = \sum\limits_{j=1}^{\mathbb{I}_{1}^{N}} \mu_j \mathcal{V}_j  \star_1 \mathcal{V}^{H}_j$. Besides, we are given a polygamma funtion $\omega^{(k)}(x)$  for any $k \in \mathbb{N}$ and $x > 0$.  For any $\theta > 0$, we have
\begin{eqnarray}\label{eq1:cor:TB Lipschitz Estimates:polygamma quasi}
\mathrm{Pr}\left(  \left\Vert \mathcal{D} \star_N \omega^{(k)}(\mathcal{A}) -  \omega^{(k)}(\mathcal{B}) \star_N \mathcal{D} \right\Vert_{\rho}  \geq \theta \right) \leq \frac{  \left( \mathbb{I}_1^N \right)^2 \Omega_{\omega^{[k+1]}  }}{\theta} \mathbb{E}\left( \left\Vert \mathcal{D} \star_N \mathcal{A} - \mathcal{B} \star_N \mathcal{D}  \right\Vert_{\rho} \right),
\end{eqnarray} 
where 
\begin{eqnarray}
\Omega_{\omega^{[k+1]}}  & = &  \max\left\{ \omega^{(k+1)}\left(\frac{1}{e}\left( \frac{(b^*)^{b^*}}{(a^*)^{a^*}}\right)^{\frac{1}{b^*-a^*}} \right)   ,  \omega^{(k+1)}(x^*) \right\}.
\end{eqnarray}
The values $a^*$ and $b^*$ are the maximizers of the function $\omega^{(k+1)}\left(\frac{1}{e}\left( \frac{b^{b}}{a^{a}}\right)^{\frac{1}{b-a}} \right)$ given $a \neq b$ and $a, b > 0$. The value $x^*$ is the maximizer of the function $\omega^{(k+1)}(x)$ given $x > 0$. 
\end{corollary}
\textbf{Proof:}

This corrollary is proved by applying Theorem 1 from~\cite{qi2008new} to Eq.~\eqref{eq:polygamma func def} and  Theorem~\ref{thm:TB Lipschitz Estimates quasi}.
$\hfill \Box$

\section{Continuity of Random Tensor Integral}\label{sec:Continuity of Random Tensor Integral}

In this section, we will establish continuity of DTI. We need the following definition to define the convergence in mean for random tensors.

\begin{definition}\label{def:conv in mean}
We say that a sequence of random tensor $\mathcal{X}_n$ converges in the $r$-th mean towards the random tensor $\mathcal{X}$ with respect to the tensor norm $\left\Vert \cdot \right\Vert_{\rho}$, if we have
\begin{eqnarray}
\mathbb{E}\left( \left\Vert \mathcal{X}_n \right\Vert^{r}_{\rho} \right)~~~\mbox{exists,}
\end{eqnarray}
and
\begin{eqnarray}
\mathbb{E}\left( \left\Vert \mathcal{X} \right\Vert^{r}_{\rho} \right)~~~ \mbox{exists,}
\end{eqnarray}
and
\begin{eqnarray}
\lim\limits_{n \rightarrow \infty}\mathbb{E}\left( \left\Vert \mathcal{X}_n  - \mathcal{X} \right\Vert^{r}_{\rho} \right) = 0. 
\end{eqnarray}
We adopt the notatation $\mathcal{X}_n \xrightarrow[]{r} \mathcal{X}$ to represent that random tensors $\mathcal{X}_n$ converges in the $r$-th mean to the random tensor $\mathcal{X}$ with respect to the tensor norm $\left\Vert \cdot \right\Vert_{\rho}$.
\end{definition}

Besides random tensor convergence definition, we also need to define triple tensor integral and second-order divide difference.

We define \emph{triple tensor integrals} (TTI) with respect to Hermitian tensors $\mathcal{A}, \mathcal{B}, \mathcal{C} \in \mathbb{C}^{I_1 \times \cdots \times I_N \times I_1 \times \cdots \times I_N}$ such that $\mathcal{A} = \sum\limits_{i=1}^{\mathbb{I}_{1}^{N}} \lambda_i \mathcal{U}_i \star_1 \mathcal{U}^{H}_i \define  \sum\limits_{i=1}^{\mathbb{I}_{1}^{N}} \lambda_i   \mathcal{P}_{\mathcal{U}_i} $, $\mathcal{B} = \sum\limits_{j=1}^{\mathbb{I}_{1}^{N}} \mu_j \mathcal{V}_j  \star_1 \mathcal{V}^{H}_j \define  \sum\limits_{j=1}^{\mathbb{I}_{1}^{N}} \mu_j   \mathcal{P}_{\mathcal{V}_j} $ and $\mathcal{C} = \sum\limits_{k=1}^{\mathbb{I}_{1}^{N}} \nu_k \mathcal{W}_j  \star_1 \mathcal{W}^{H}_j \define  \sum\limits_{k=1}^{\mathbb{I}_{1}^{N}} \nu_k   \mathcal{P}_{\mathcal{W}_k} $. Given the function $\varphi: \mathbb{R}^3 \rightarrow \mathbb{R}$, the TTI associated with tensors $\mathcal{A}, \mathcal{B}, \mathcal{C} $ and the function $\varphi$, denoted as $T_{\mathcal{A}, \mathcal{B}, \mathcal{C}, \varphi}( \mathcal{X}, \mathcal{Y} )$, can be expressed as
\begin{eqnarray}\label{eq:Triple Tensor Int Def}
T_{\mathcal{A}, \mathcal{B}, \mathcal{C}, \varphi}(\mathcal{X}, \mathcal{Y}) &=& \sum\limits_{i=1}^{\mathbb{I}_{1}^{N}} \sum\limits_{j=1}^{\mathbb{I}_{1}^{N}} \sum\limits_{k=1}^{\mathbb{I}_{1}^{N}} \varphi( \lambda_i, \mu_j, \nu_k)  \mathcal{P}_{\mathcal{U}_i} \star_N \mathcal{X} \star_N \mathcal{P}_{\mathcal{V}_j} \star_N \mathcal{Y} \star_N \mathcal{P}_{\mathcal{W}_j},
\end{eqnarray}
where $\mathcal{X}, \mathcal{Y} \in \mathbb{C}^{I_1 \times \cdots \times I_N \times I_1 \times \cdots \times I_N}$. 

The second-order divide difference for a function $f(x)$, denoted as $f^{[2]}(x,y,z)$, can be defined as
\begin{eqnarray}\label{eq:2 order divide diff def}
f^{[2]}(x,y,z) \define \frac{f^{[1]}(y,z) -  f^{[1]}(x,y) }{z - x}
\end{eqnarray}

\begin{lemma}\label{lma:bound for Triple Tensor Int}
Given three Hermitian tensors $\mathcal{A}, \mathcal{B}, \mathcal{C} \in \mathbb{C}^{I_1 \times \cdots \times I_N \times I_1 \times \cdots \times I_N}$ such that $\mathcal{A} = \sum\limits_{i=1}^{\mathbb{I}_{1}^{N}} \lambda_i \mathcal{U}_i \star_1 \mathcal{U}^{H}_i \define  \sum\limits_{i=1}^{\mathbb{I}_{1}^{N}} \lambda_i   \mathcal{P}_{\mathcal{U}_i} $, $\mathcal{B} = \sum\limits_{j=1}^{\mathbb{I}_{1}^{N}} \mu_j \mathcal{V}_j  \star_1 \mathcal{V}^{H}_j \define  \sum\limits_{j=1}^{\mathbb{I}_{1}^{N}} \mu_j   \mathcal{P}_{\mathcal{V}_j} $ and $\mathcal{C} = \sum\limits_{k=1}^{\mathbb{I}_{1}^{N}} \nu_k \mathcal{W}_j  \star_1 \mathcal{W}^{H}_j \define  \sum\limits_{k=1}^{\mathbb{I}_{1}^{N}} \nu_k   \mathcal{P}_{\mathcal{W}_k} $ and the function $f(x)$ with $f^{[2]}(x)$ bounded by $\Omega_{f^{[2]}}$, we then have the following norm estimate for $T_{\mathcal{A}, \mathcal{B}, \mathcal{C}, \varphi}( \mathcal{X}, \mathcal{Y} )$:
\begin{eqnarray}
\left\Vert T_{\mathcal{A}, \mathcal{B}, \mathcal{C}, \varphi}(\mathcal{X}, \mathcal{Y}) \right\Vert_{\rho}   &\leq& 
\left( \mathbb{I}_{1}^{N}  \right)^3  \Omega_{f^{[2]}}  \left\Vert \mathcal{X} \right\Vert_{\rho} \cdot \left\Vert \mathcal{Y} \right\Vert_{\rho}. 
\end{eqnarray}
\end{lemma}
\textbf{Proof:}

Since we have
\begin{eqnarray}
\left\Vert T_{\mathcal{A}, \mathcal{B}, \mathcal{C}, \varphi}(\mathcal{X}, \mathcal{Y}) \right\Vert_{\rho} &=& \left\Vert 
\sum\limits_{i=1}^{\mathbb{I}_{1}^{N}} \sum\limits_{j=1}^{\mathbb{I}_{1}^{N}} \sum\limits_{k=1}^{\mathbb{I}_{1}^{N}} \varphi( \lambda_i, \mu_j, \nu_k)  \mathcal{P}_{\mathcal{U}_i} \star_N \mathcal{X} \star_N \mathcal{P}_{\mathcal{V}_j} \star_N \mathcal{Y} \star_N \mathcal{P}_{\mathcal{W}_j}  \right\Vert_{\rho} \nonumber \\
& \leq_1 &
\sum\limits_{i=1}^{\mathbb{I}_{1}^{N}} \sum\limits_{j=1}^{\mathbb{I}_{1}^{N}} \sum\limits_{k=1}^{\mathbb{I}_{1}^{N}}
 \left\vert \varphi( \lambda_i, \mu_j, \nu_k)  \right\vert  \left\Vert   \mathcal{P}_{\mathcal{U}_i} \star_N \mathcal{X} \star_N \mathcal{P}_{\mathcal{V}_j} \star_N \mathcal{Y} \star_N \mathcal{P}_{\mathcal{W}_j} \right\Vert_{\rho} \nonumber \\
& \leq_2 & \sum\limits_{i=1}^{\mathbb{I}_{1}^{N}} \sum\limits_{j=1}^{\mathbb{I}_{1}^{N}}   \sum\limits_{k=1}^{\mathbb{I}_{1}^{N}}
 \left\vert \varphi( \lambda_i, \mu_j, \nu_k)  \right\vert   \left\Vert \mathcal{X} \right\Vert_{\rho} \cdot \left\Vert \mathcal{Y} \right\Vert_{\rho} \nonumber \\
& \leq & \left( \mathbb{I}_{1}^{N}  \right)^3  \Omega_{f^{[2]}}  \left\Vert \mathcal{X} \right\Vert_{\rho} \cdot \left\Vert \mathcal{Y} \right\Vert_{\rho} 
\end{eqnarray}
where $\leq_1$ comes from triangle inequality of the unitarily invariant norm and $\leq_2$ comes from the definition of the unitarily invariant norm and submultiplicative property of any unitarily invariant norm.  
$\hfill \Box$

Following theorem is about the continuity of a tensor integral. 

\begin{theorem}\label{thm:continuity of double TI}
Let $\mathcal{A}_n, \mathcal{A}, \mathcal{B}_n, \mathcal{B} \in \mathbb{C}^{I_1 \times \cdots \times I_N \times I_1 \times \cdots \times I_N}$ be random Hermitian tensors such that
\begin{eqnarray}\label{eq1:thm:continuity of double TI}
\mathcal{A}_n \xrightarrow[]{r} \mathcal{A},~\mbox{and}~
\mathcal{B}_n \xrightarrow[]{r} \mathcal{B},
\end{eqnarray}
where $1 \leq r < \infty$. Moreover, a real values function $f(x)$ for $x \in \mathbb{R}$ such that $ f^{[2]}(x)$ exists and bounded by $\Omega_{f^{[2]}}$, respectively. Then, we have  
\begin{eqnarray}
T_{\mathcal{A}_n, \mathcal{B}_n, f^{[1]}}(\mathcal{X}) \xrightarrow[]{r}
T_{\mathcal{A}, \mathcal{B}, f^{[1]}}(\mathcal{X}),
\end{eqnarray}
where $\mathcal{X} \in \mathbb{C}^{I_1 \times \cdots \times I_N \times I_1 \times \cdots \times I_N}$ is a fixed tensor. 
\end{theorem}
\textbf{Proof:}

From Lemma~\ref{lma:pert formula} and telescoping summation, we have the following:
\begin{eqnarray}\label{eq2:thm:continuity of double TI}
\left\Vert T_{\mathcal{A}_n, \mathcal{B}_n, f^{[1]}}(\mathcal{X}) - T_{\mathcal{A}, \mathcal{B}, f^{[1]}}(\mathcal{X})\right\Vert_{\rho} &=& \left\Vert T_{\mathcal{A}_n, \mathcal{B}_n, f^{[1]}}(\mathcal{X}) - T_{\mathcal{A}, \mathcal{B}_n, f^{[1]}}(\mathcal{X})  \right. \nonumber \\
&    & \left. + T_{\mathcal{A}, \mathcal{B}_n, f^{[1]}}(\mathcal{X}) -  T_{\mathcal{A}, \mathcal{B}, f^{[1]}}(\mathcal{X}) \right\Vert_{\rho} \nonumber \\
&=_1& \left\Vert T_{\mathcal{A}_n, \mathcal{A}, \mathcal{B}_n, f^{[2]}}(\mathcal{A}_n - \mathcal{A}, \mathcal{X}) + T_{\mathcal{A}, \mathcal{B}_n, \mathcal{B}, f^{[2]}}(X, \mathcal{B}_n - \mathcal{B}) \right\Vert_{\rho} \nonumber \\
&\leq& \Omega_{f^{[2]}} \left\Vert \mathcal{A}_n - \mathcal{A} \right\Vert_{\rho}\left\Vert  \mathcal{X} \right\Vert_{\rho}  + \Omega_{f^{[2]}} \left\Vert \mathcal{B}_n - \mathcal{B} \right\Vert_{\rho}\left\Vert  \mathcal{X} \right\Vert_{\rho},
\end{eqnarray}
where $=_1$ is obtained from the $f^{[2]}$ definition given by Eq.~\eqref{eq:2 order divide diff def}, and the $\leq$ comes from Lemma~\ref{lma:bound for Triple Tensor Int} and triangle inequality.

By raising the power $r$ and taking the expectation at the both sides of the inequality provided by Eq.~\eqref{eq2:thm:continuity of double TI}, we have proved this theorem by conditions given by Eq.~\eqref{eq1:thm:continuity of double TI} and the following inequality:
\begin{eqnarray}
(a+b)^r \leq 2^r(a^r + b^r)~~\mbox{given $a, b \geq 0$}.
\end{eqnarray}
$\hfill \Box$

\section{Applications of Tensor Integral}\label{sec:Applications of Tensor Integral}

In this section, we will apply Theorem~\ref{thm:continuity of double TI} and perturbation formulas provided by Lemma~\ref{lma:pert formula} and Lemma~\ref{lma:pert formula quasi} to bound the tail probability of the derivative of tensor-valued function norm. 

Given a fixed perturbation tensor $\mathcal{X} \in \mathbb{C}^{I_1 \times \cdots \times I_N \times I_1 \times \cdots \times I_N}$ with respect to the random Hermitian tensor  $\mathcal{A} \in \mathbb{C}^{I_1 \times \cdots \times I_N \times I_1 \times \cdots \times I_N}$, and a tensor-valued function $f(x)$, we define the derivative of $f(x)$ at $\mathcal{A}$ with respect to the perturbation $\mathcal{X}$, represented by $f'_{\mathcal{X}}(\mathcal{A})$, as
\begin{eqnarray}\label{eq:derivative def}
f'_{\mathcal{X}}(\mathcal{A}) \define \lim\limits_{t \rightarrow 0} \frac{f(  \mathcal{A} + t \mathcal{X}) - f(\mathcal{A}) }{t }
\end{eqnarray}

Following theorem is about the tail bound for the norm of $f'_{\mathcal{X}}(\mathcal{A})$.

\begin{theorem}\label{thm:derivative TB}
Given a fixed perturbation tensor $\mathcal{X} \in \mathbb{C}^{I_1 \times \cdots \times I_N \times I_1 \times \cdots \times I_N}$ and a random Hermitian tensor  $\mathcal{A} \in \mathbb{C}^{I_1 \times \cdots \times I_N \times I_1 \times \cdots \times I_N}$ with a tensor-valued function $f(x)$. Suppose we have $\mathcal{A} = \sum\limits_{i=1}^{\mathbb{I}_{1}^{N}} \lambda_i \mathcal{P}_{\mathcal{U}_i}$.Then, we have the tail bound for $\left\Vert f'_{\mathcal{X}}(\mathcal{A}) \right\Vert_{\rho}$ as
\begin{eqnarray}\label{eq1:thm:derivative TB}
\mathrm{Pr}\left( \left\Vert f'_{\mathcal{X}}(\mathcal{A}) \right\Vert_{\rho} \geq \theta \right) &\leq& 
\frac{  \left( \mathbb{I}_1^N \right)^2 \left\Vert \mathcal{X} \right\Vert_{\rho} }{\theta}  \sum\limits_{i=1}^{\mathbb{I}_{1}^{N}}  \sum\limits_{j=1}^{\mathbb{I}_{1}^{N}} \mathbb{E}\left( \left\vert  f^{[1]}(\lambda_i, \lambda_j)\right\vert \right).
\end{eqnarray}
\end{theorem}
\textbf{Proof:}

From Lemma~\ref{lma:pert formula}, we have
\begin{eqnarray}
f(  \mathcal{A} + t \mathcal{X}) - f(\mathcal{A})  &=& T_{\mathcal{A}+  t \mathcal{X}, \mathcal{A}, f^{[1]}}\left( \mathcal{A} +  t \mathcal{X} - \mathcal{A} \right).  
\end{eqnarray}
Then, we have
\begin{eqnarray}
\frac{f(  \mathcal{A} + t \mathcal{X}) - f(\mathcal{A})}{t}  &=& T_{\mathcal{A}+  t \mathcal{X}, \mathcal{A}, f^{[1]}}\left( \mathcal{X} \right).  
\end{eqnarray}

From Theorem~\ref{thm:continuity of double TI}, we have
\begin{eqnarray}
f'_{\mathcal{X}}(\mathcal{A}) &=& \lim\limits_{t \rightarrow 0} \frac{f(  \mathcal{A} + t \mathcal{X}) - f(\mathcal{A}) }{t } \nonumber \\
& \xrightarrow[]{r} & T_{\mathcal{A}, \mathcal{A}, f^{[1]}}\left( \mathcal{X} \right).  
\end{eqnarray}
Since the convergence in the $r$-th mean implies the convergence in probability, we then have
\begin{eqnarray}\label{eq2:thm:derivative TB}
\mathrm{Pr}\left( \left\Vert f'_{\mathcal{X}}(\mathcal{A}) \right\Vert_{\rho} \geq \theta \right) &=& 
\mathrm{Pr}\left( \left\Vert T_{\mathcal{A}, \mathcal{A}, f^{[1]}}\left( \mathcal{X} \right) \right\Vert_{\rho} \geq \theta \right) \nonumber \\
&\leq&  \frac{  \left( \mathbb{I}_1^N \right)^2 \left\Vert \mathcal{X} \right\Vert_{\rho} }{\theta}  \sum\limits_{i=1}^{\mathbb{I}_{1}^{N}}  \sum\limits_{j=1}^{\mathbb{I}_{1}^{N}} \mathbb{E}\left( \left\vert  f^{[1]}(\lambda_i, \lambda_j)\right\vert \right), 
\end{eqnarray}
where the inequality comes from Theorem~\ref{thm:unitary inv norm tail bound}
$\hfill \Box$

The commutator of a tensor $\mathcal{A} \in \mathbb{C}^{I_1 \times \cdots \times I_N \times I_1 \times \cdots \times I_N} $ with respect to the tensor $\mathcal{B} \in \mathbb{C}^{I_1 \times \cdots \times I_N \times I_1 \times \cdots \times I_N} $ if we have $\mathcal{A} \star_N \mathcal{B} = \mathcal{B} \star_N \mathcal{A}$.

Given a fixed perturbation tensor $\mathcal{X} \in \mathbb{C}^{I_1 \times \cdots \times I_N \times I_1 \times \cdots \times I_N}$ with respect to the random Hermitian tensor  $\mathcal{A} \in \mathbb{C}^{I_1 \times \cdots \times I_N \times I_1 \times \cdots \times I_N}$, and a tensor-valued function $f(x)$, we define the derivative of $f(x)$ at $\mathcal{A}$ with respect to the perturbation $\mathcal{X}$ and the commutator tensor $\mathcal{D} \in \mathbb{C}^{I_1 \times \cdots \times I_N \times I_1 \times \cdots \times I_N}$ of $\mathcal{A}$, represented by $f'_{\mathcal{X} | \mathcal{D} }(\mathcal{A})$, as
\begin{eqnarray}\label{eq:derivative def quasi}
f'_{\mathcal{X}| \mathcal{D}}(\mathcal{A}) \define \lim\limits_{t \rightarrow 0} \frac{\mathcal{D} \star_N f(  \mathcal{A} + t \mathcal{X}) - f(\mathcal{A}) \star_N \mathcal{D}}{t }
\end{eqnarray}

Following theorem is about the tail bound for the norm of $f'_{\mathcal{X}| \mathcal{D} }(\mathcal{A})$.

\begin{theorem}\label{thm:derivative TB quasi}
Given a fixed perturbation tensor $\mathcal{X} \in \mathbb{C}^{I_1 \times \cdots \times I_N \times I_1 \times \cdots \times I_N}$ and a random Hermitian tensor  $\mathcal{A} \in \mathbb{C}^{I_1 \times \cdots \times I_N \times I_1 \times \cdots \times I_N}$ with a tensor-valued function $f(x)$. Suppose we have $\mathcal{A} = \sum\limits_{i=1}^{\mathbb{I}_{1}^{N}} \lambda_i \mathcal{P}_{\mathcal{U}_i}$ and the tensor $\mathcal{D} \in \mathbb{C}^{I_1 \times \cdots \times I_N \times I_1 \times \cdots \times I_N}$ is the commutator of the tensor $\mathcal{A}$. Then, we have the tail bound for $\left\Vert f'_{\mathcal{X} | \mathcal{D} }(\mathcal{A}) \right\Vert_{\rho}$ as
\begin{eqnarray}\label{eq1:thm:derivative TB quasi}
\mathrm{Pr}\left( \left\Vert f'_{\mathcal{X} | \mathcal{D} }(\mathcal{A}) \right\Vert_{\rho} \geq \theta \right) &\leq& 
\frac{  \left( \mathbb{I}_1^N \right)^2 \left\Vert \mathcal{D} \star_N \mathcal{X} \right\Vert_{\rho} }{\theta}  \sum\limits_{i=1}^{\mathbb{I}_{1}^{N}}  \sum\limits_{j=1}^{\mathbb{I}_{1}^{N}} \mathbb{E}\left( \left\vert  f^{[1]}(\lambda_i, \lambda_j)\right\vert \right).
\end{eqnarray}
\end{theorem}
\textbf{Proof:}

From Lemma~\ref{lma:pert formula quasi}, we have
\begin{eqnarray}
\mathcal{D}\star_N f(  \mathcal{A} + t \mathcal{X}) - f(\mathcal{A}) \star_N \mathcal{D}  &=& T_{\mathcal{A}+  t \mathcal{X}, \mathcal{A}, f^{[1]}}\left( \mathcal{D} \star_N \left( \mathcal{A} +  t \mathcal{X} \right) - \mathcal{A} \star_N \mathcal{D} \right).  
\end{eqnarray}
Also, from $\mathcal{D}\star_N  \mathcal{A} = \mathcal{A} \star_N \mathcal{D}$, we also have 
\begin{eqnarray}
\frac{ \mathcal{D}\star_N f(  \mathcal{A} + t \mathcal{X}) - f(\mathcal{A})\star_N \mathcal{D}   }{t}  &=& T_{\mathcal{A}+  t \mathcal{X}, \mathcal{A}, f^{[1]}}\left( \mathcal{D} \star_N \mathcal{X} \right).  
\end{eqnarray}

From Theorem~\ref{thm:continuity of double TI}, we have
\begin{eqnarray}
f'_{\mathcal{X} | \mathcal{D}}(\mathcal{A}) &=& \lim\limits_{t \rightarrow 0} \frac{  \mathcal{D} \star_N f( \mathcal{A} + t \mathcal{X}) - f(\mathcal{A}) \star_N \mathcal{D} }{t } \nonumber \\
& \xrightarrow[]{r} & T_{\mathcal{A}, \mathcal{A}, f^{[1]}}\left( \mathcal{D} \star_N \mathcal{X} \right).  
\end{eqnarray}
Since the convergence in the $r$-th mean implies the convergence in probability, we then have
\begin{eqnarray}\label{eq2:thm:derivative TB quasi}
\mathrm{Pr}\left( \left\Vert f'_{\mathcal{X} | \mathcal{D} }(\mathcal{A}) \right\Vert_{\rho} \geq \theta \right) &=& 
\mathrm{Pr}\left( \left\Vert T_{\mathcal{A}, \mathcal{A}, f^{[1]}}\left( \mathcal{D} \star_N \mathcal{X} \right) \right\Vert_{\rho} \geq \theta \right) \nonumber \\
&\leq&  \frac{  \left( \mathbb{I}_1^N \right)^2 \left\Vert \mathcal{D} \star_N \mathcal{X} \right\Vert_{\rho} }{\theta}  \sum\limits_{i=1}^{\mathbb{I}_{1}^{N}}  \sum\limits_{j=1}^{\mathbb{I}_{1}^{N}} \mathbb{E}\left( \left\vert  f^{[1]}(\lambda_i, \lambda_j)\right\vert \right), 
\end{eqnarray}
where the inequality comes from Theorem~\ref{thm:unitary inv norm tail bound} again. 
$\hfill \Box$

\section{Conclusions}\label{sec:Conclusions}

We first define what is the random DTI and derive the tail bound of the unitarily invariant norm for a random DTI. This bound assists us to establish tail bounds of the unitarily invariant norm for various types of dual tensor means, e.g., arithmetic mean, geometric mean, harmonic mean, and general mean. The random DTI is also being applied to build the random Lipschitz estimate in contexts of random tensors. Finally, we derive the continuity property for random DTI in the sense of convergence in the random tensor mean, and apply this fact to obtain the tail bound of the unitarily invariant norm for the derivative of the tensor-valued function. Possible future works will be to extend DTI to multiple tensor integrals. 

\bibliographystyle{IEEETran}
\bibliography{PDTensor_Int_Bib}

\end{document}